\documentclass[a4paper,10pt]{article}
\usepackage{stmaryrd}
\usepackage{amsfonts}
\usepackage{bbm}
\usepackage{amscd}
\usepackage{mathrsfs}
\usepackage{latexsym,amssymb,amsmath,amscd,amscd,amsthm,amsxtra,xypic}
\usepackage[dvips]{graphicx}
\usepackage[utf8]{inputenc}
\usepackage[T1]{fontenc}
\usepackage{lmodern}
\usepackage{amssymb}
\usepackage[all]{xy}
\usepackage{nicefrac,mathtools,enumitem}
\usepackage{microtype}

\newtheorem{thm}{Theorem}[section]
\newtheorem{lem}[thm]{Lemma}
\newtheorem{cor}[thm]{Corollary}
\newtheorem{pro}[thm]{Proposition}

\newtheorem{rmk}[thm]{Remark}
\newtheorem{defi}[thm]{Definition}

\newcommand{\be }{\begin{equation}}
\newcommand{\ee }{\end{equation}}

\newcommand{\pf}{\noindent{\bf Proof.}\ }

\newcommand{\CWM}{C^{\infty}(M)}

\newcommand{\frkG}{\mathfrak G}

\newcommand{\frkX}{\mathfrak X}

\def\qed{\hfill ~\vrule height6pt width6pt depth0pt}

\newcommand{\half}{\frac{1}{2}}

\newcommand{\pairing}[1]{\left ( #1\right )}

\newcommand{\Courant}[1]{\left\llbracket  #1\right\rrbracket }
\newcommand{\Dorfman}[1]{\{ #1\}}

\newcommand{\dev}{\mathfrak{D}}

\newcommand{\dM}{d}

\newcommand{\Ham}{\mathrm{Ham}}

\begin{document}
\title{
{On higher analogues of Courant algebroids
\thanks
 {
Research partially supported by NSF of China (10871007), US-China
CMR Noncommutative Geometry (10911120391/A0109), China Postdoctoral
Science Foundation (20090451267),  Science Research Foundation for
Excellent Young Teachers of Mathematics School at Jilin University.
 }
} }
\author{
Yanhui Bi \\ Department of Mathematics and LMAM,\\\vspace{2mm}
Peking University, Beijing 100871, China\\
Yunhe Sheng  \\
Mathematics School $\&$ Institute of Jilin University,\\
 Changchun 130012, Jilin, China
\\AND\\
 Department of Mathematics, Dalian University of
Technology,\\
 Dalian 116023,  China\\
email: shengyh@jlu.edu.cn}

\date{}


\maketitle
\begin{abstract}
In this paper, we study the algebraic properties of the higher
analogues of Courant algebroid structures on the direct sum bundle
$TM\oplus\wedge^nT^*M$ for an $m$-dimensional manifold. As an
application,  we revisit Nambu-Poisson structures and
multisymplectic structures. We prove that the graph of an
$(n+1)$-vector field $\pi$ is closed under the higher-order Dorfman
bracket iff  $\pi$ is a Nambu-Poisson structure. Consequently, there
is an induced Leibniz algebroid structure on $\wedge^nT^*M$. The
graph of an $(n+1)$-form $\omega$ is closed under the higher-order
Dorfman bracket iff $\omega$ is a premultisymplectic structure of
order $n$, i.e. $\dM\omega=0$. Furthermore, there is a Lie algebroid
structure on the admissible bundle $A\subset\wedge^{n}T^*M$. In
particular, for a $2$-plectic structure, it induces the Lie
2-algebra structure given in \cite{baez:classicalstring}.

\end{abstract}
\section{Introduction}

The notion of Courant algebroid was introduced in \cite{lwx}  to
study the double of Lie bialgebroids. Equivalent definition was
given by Roytenberg  \cite{Roytenbergphdthesis}. In resent years,
Courant algebroids have been far and wide studied from several
aspects and have been found
 many applications in the theory of Manin pairs and moment
maps \cite{alekseevmanin,Liblandmeinrenken}; generalized complex
structures \cite{gualtieri,xureductionofgeneralized};
$L_{\infty}$-algebras and symplectic supermanifolds
\cite{Roytenbergphdthesis}; gerbes
 as well as BV algebras and  topological field
theories.

Roughly speaking, a Courant algebroid is a vector bundle together
with a nondegenerate symmetric pairing, a bracket operation and an
anchor map such that some compatible conditions are satisfied. The
standard Courant algebroid is the direct sum bundle $TM\oplus T^*M$.
The standard Courant bracket is given by
$$
\Courant{X+\alpha,Y+\beta}=[X,Y]+L_X\beta-L_Y\alpha+\half(\dM
i_Y\alpha-\dM i_X\beta).
$$
However, many experts know that on the direct sum bundle $TM\oplus
\wedge^nT^*M$, there is also a similar bracket operation, which we
call the  higher-order Courant bracket. But the nature of its full
algebraic structures
 remains a work in progress. The main purpose of this paper is to study
 its properties. We proved that the Jacobi identity holds up to an
 exact term. Moreover, we can also introduce the Dorfman bracket
 $\Dorfman{\cdot,\cdot}$ which is not skew-symmetric. It follows that $(TM\oplus
 \wedge^nT^*M,\Dorfman{\cdot,\cdot},\rho)$ is a Leibniz algebroid,
 where $\rho $ is the projection to $TM$. As an application,  we revisit Nambu-Poisson structures and multisymplectic
 structures.

 Nambu-Poisson structures were
introduced in \cite{TakhtajanNambu} by Takhtajan in order to find an
axiomatic formalism for Nambu-mechanics which is a generalization of
Hamiltonian mechanics \cite{chatterjee1,chatterjee2}. The local and
global structures of a Nambu-Poisson manifold were studied in
\cite{MarrerodynamicsNambu,MarmonPoisson,NakanishiNM}. It was
further studied by Grabowski and Marmo \cite{GrabowskiNambu2}. In
\cite{MarreroLeibnizNambu}, the authors proved that associated to a
Nambu-Poisson structure\footnote{it is an $(n+1)$-vector field
satisfying  some integrability conditions} of order $n$, their is a
bracket operation  on $\wedge^nT^*M$ such that it is a Leibniz
algebroid. However, we do not think it is the best bracket that make
$\wedge^nT^*M$ to be a Leibniz algebroid. We find
 that $\pi$ is a Nambu-Poisson structure of order $n$ if and only if
 the graph of $\pi$ is closed under the higher-order Dorfman bracket operation
 $\Dorfman{\cdot,\cdot}$. Consequently, there is an induced  bracket
operation  $[\cdot,\cdot]_\pi$ on $\wedge^nT^*M$,
$$
[\alpha,\beta]_\pi=L_{\pi^\sharp(\alpha)}\beta-L_{\pi^\sharp(\beta)}\alpha+\dM
i_{\pi^\sharp(\beta)}\alpha,
$$
such that $(\wedge^nT^*M,[\cdot,\cdot]_\pi,\rho)$ is a Leibniz
algebroid, where $\rho$ is the projection to $TM$. In particular, if
$n=1$, we recover the cotangent bundle Lie algebroid.

In \cite{baez:classicalstring}, Baez, Hoffnung and Rogers introduced
a hemi-bracket and a semi-bracket on the set of Hamiltonian
$(n-1)$-forms associated to a multisymplectic structure of order
$n$, or equivalently, an $n$-plectic structure\footnote{it is an
$(n+1)$-form which is closed and nondegenerate}. In particular, one
can obtain a Lie 2-algebra for a 2-plectic structure. Recently,
several aspects of the theory have been studied in depth
\cite{Rogers,zambonCourant}. In this paper, we introduce a Lie
bracket on admissible $n$-forms. Since admissible $n$-forms are
sections  of a vector bundle $A\subset \wedge^n T^*M$, we obtain a
Lie algebroid. We also show that how one can obtain Baez's
semi-bracket and hemi-bracket on Hamiltonian $(n-1)$-forms from our
Lie algebroid. The relation between $n$-plectic structures and the
Leibniz algebroid $(TM\oplus
 \wedge^nT^*M,\Dorfman{\cdot,\cdot},\rho)$ is also discussed.

The paper is organized as follows. In Section 2  we study the
algebraic properties of the higher-order Courant bracket
$\Courant{\cdot,\cdot}$ and the higher-order Dorfman bracket
$\Dorfman{\cdot,\cdot}$ on the direct sum bundle $TM\oplus
\wedge^nT^*M$. We prove that the Jacobi identity of
$\Courant{\cdot,\cdot}$ holds up to an exact term and the triple
$(TM\oplus
 \wedge^nT^*M,\Dorfman{\cdot,\cdot},\rho)$ is a Leibniz
algebroid. In Section 3 we study the relation between Nambu-Poisson
structures and the higher-order Dorfman bracket. There is an induced
Leibniz algebroid structure on $\wedge^nT^*M$ associated to a
Nambu-Poisson structure.  In Section 4 we revisit multisymplectic
structures and construct a Lie algebroid, whose sections are
admissible $n$-forms.\vspace{3mm}

{\bf Acknowledgement:} We give warmest thanks to Zhangju Liu and
Chenchang Zhu for useful discussion. We also thank the referees for
very helpful comments.

\section{Higher-order Courant brackets}
In this section, we consider the higher-order Courant brackets on
the direct sum bundle $$\mathcal{T}^n\triangleq
TM\oplus\wedge^nT^*M$$ and associated properties. There is a natural
$\wedge^{n-1}T^*M$-valued nondegenerate symmetric pairing
$\pairing{\cdot,\cdot}$ on it,
\begin{equation}\label{pair}
\pairing{X+\alpha,Y+\beta}=\half(i_X\beta+i_Y\alpha),\quad\forall~X,Y\in\frkX(M),~\alpha,\beta\in\Omega^n(M).
\end{equation}
The higher-order Courant bracket  $\Courant{\cdot,\cdot}$  on the
section space of $\mathcal{T}^n$  is defined as follows,
\begin{equation}\label{Cbracket}
\Courant{X+\alpha,Y+\beta}=[X,Y]+L_X\beta-L_Y\alpha+\half(\dM
i_Y\alpha-\dM i_X\beta).
\end{equation}
Obviously, it is skew-symmetric. When $n=0,$ we obtain the Lie
algebroid $TM\times \mathbb R$. When $n=1$, it is exactly the
standard Courant bracket \cite{lwx,Roytenbergphdthesis}. Let
$\rho:\mathcal{T}^n\longrightarrow TM$ be the projection, i.e.
\begin{equation}\label{anchor}
\rho(X+\alpha)=X,\quad\forall ~X+\alpha\in \Gamma(\mathcal{T}^n).
\end{equation}

We call
$(\mathcal{T}^n,\pairing{\cdot,\cdot},\Courant{\cdot,\cdot},\rho)$
the higher analogues of Courant algebroids.
\begin{rmk}
 In \cite{rw}, Roytenberg and
Weinstein proved that a Courant algebroid gives rise to a $2$-term
$L_\infty$-algebra. However, as proposed in Section 3.8 in
\cite{gualtieri}, the nature of the full algebraic structures for
all $TM\oplus\wedge^nT^*M$ remains a work in progress.  For more
details about generalized complex structures of Courant-Jacobi
algebroids, see \cite{shengJacobi}
\end{rmk}

The higher-order Courant bracket satisfies some similar properties
as the Courant bracket.
\begin{thm}
For any $e_1,e_2,e_3\in \Gamma(\mathcal{T}^n),~f\in
\CWM,~\xi\in\Omega^{n-1}(M)$,  we have
\begin{itemize}
\item[\rm(i).]
$ \Courant{e_1,\Courant{e_2,e_3}}+c.p.=\dM T(e_1,e_2,e_3), $ where
$T:\wedge^3\mathcal{T}^n\longrightarrow\wedge^{n-1}T^*M$ is defined
by
$$
T(e_1,e_2,e_3)=-\frac{1}{3}\big(\langle\Courant{e_1,e_2},e_3\rangle+c.p.\big).
$$
\item[\rm(ii).]$\Courant{e_1,fe_2}=f\Courant{e_1,e_2}+\rho(e_1)(f)e_2-\dM f\wedge(e_1,e_2).$
\item[\rm(iii).]$\rho\Courant{e_1,e_2}=[\rho(e_1),\rho(e_2)].$
\item[\rm(iv).]$L_{\rho(e_1)}(e_2,e_3)=(\Courant{e_1,e_2}+d(e_1,e_2),e_3)+(e_2,\Courant{e_1,e_3}+d(e_1,e_3)).$
\end{itemize}
\end{thm}

 The properties $(i)-(iii)$ are the same as the ones for
Courant algebroids. However, the property $(iv)$ is different since
the Lie derivation is involved in.

\begin{rmk}
  The proof of this theorem needs some computation. In the
  following, we will introduce the higher order Dorfman bracket (see
  \eqref{CCD}), which will be used more often since it satisfies
  better properties (Theorem \ref{thm:l}). To avoid repeat, we only
  give the proof of Theorem \ref{thm:l}.
\end{rmk}

 We can also introduce the following higher-order Dorfman bracket,
\begin{eqnarray}\label{CCD}
\Dorfman{e_1,e_2}=\Courant{e_1,e_2}+d(e_1,e_2),\quad\forall~e_1,e_2\in
\Gamma(\mathcal{T}^n).
\end{eqnarray}
When $n=1$, it is exactly the Dorfman bracket. Assume that
$e_1=X+\alpha,~e_2=Y+\beta$, we have
\begin{eqnarray}\label{CDbracket}
\Dorfman{X+\alpha,Y+\beta}= [X,Y]+L_X\beta-L_Y\alpha+\dM i_Y\alpha.
\end{eqnarray}
It is obvious that
$$
\Dorfman{X+\alpha,X+\alpha}=d(X+\alpha,X+\alpha)=\dM i_X\alpha.
$$

Leibniz algebroids are generalizations of Lie algebroids.

\begin{defi}{\rm\cite{MarreroLeibnizNambu}}
A Leibniz algebroid structure on a vector bundle $E\longrightarrow
M$ is a pair that consists of a Leibniz algebra structure
$\Dorfman{\cdot,\cdot}$ on the section space $\Gamma(E)$ and a
vector bundle morphism $\rho:E\longrightarrow TM$, called the
anchor, such that the following relations are satisfied:
\begin{itemize}
\item[\rm(1)]$ ~\rho\Dorfman{X,Y}=[\rho(X),\rho(Y)],\quad \forall~X,Y\in\Gamma(E)$;
\item[\rm(2)] $~\Dorfman{X,fY}=f[X,Y]+\rho(X)(f)Y,\quad \forall~f\in
\CWM.$
\end{itemize}
\end{defi}
If the bracket $\Dorfman{\cdot,\cdot}$ is skew-symmetric, we recover
the notion of Lie algebroids.\vspace{2mm}

The higher-order Dorfman bracket also satisfies similar properties
as the usual Dorfman bracket.

\begin{thm}\label{thm:l}
\begin{itemize}
\item[\rm(1).] For any $e_1,e_2\in \Gamma(\mathcal{T}^n),~f\in \CWM$, we have
\begin{eqnarray*}
\Dorfman{e_1,fe_2}&=&f\Dorfman{e_1,e_2}+\rho(e_1)(f)e_2,\\
\Dorfman{fe_1,e_2}&=&f\Dorfman{e_1,e_2}-\rho(e_2)(f)e_1+\dM f\wedge2
(e_1,e_2).
\end{eqnarray*}
\item[\rm(2).]The Dorfman bracket $\Dorfman{\cdot,\cdot}$ is a Leibniz bracket, i.e. for any $e_1,e_2,e_3\in \Gamma(\mathcal{T}^n)$,
$$
\Dorfman{e_1,\Dorfman{e_2,e_3}}=\Dorfman{\Dorfman{e_1,e_2},e_3}+\Dorfman{e_2,\Dorfman{e_1,e_3}}.
$$ Consequently, $(\mathcal{T}^n,\Dorfman{\cdot,\cdot},\rho)$ is a
Leibniz algebroid.
\item[\rm(3).]The pairing \eqref{pair} and the higher-order Dorfman bracket is
compatible in the following sense,
\begin{equation}\label{eqn:relation}
L_{\rho(e_1)}\pairing{e_2,e_3}=\pairing{\Dorfman{e_1,e_2},e_3}+\pairing{e_2,\Dorfman{e_1,e_3}}.
\end{equation}
\end{itemize}
\end{thm}
\pf (1). Write $e_1=X+\alpha,~e_2=Y+\beta$, where $X,Y\in\frkX(M)$
and $\alpha,\beta\in\Omega^1(M)$, then we have
\begin{eqnarray*}
\Dorfman{X+\alpha,f(Y+\beta)}&=&[X,fY]+L_Xf\beta-i_{fY}\dM\alpha\\
&=&f[X,Y]+X(f)Y+fL_X\beta+X(f)\beta-fi_{Y}\dM\alpha\\
&=&f\Dorfman{X+\alpha,Y+\beta}+X(f)(Y+\beta).
\end{eqnarray*}
Similarly, we can obtain the expression of $\Dorfman{fe_1,e_2}$ as
in the theorem.

 (2). Write $e_3=Z+\gamma$, where $Z\in\frkX(M)$
and $\gamma\in\Omega^1(M)$, by straightforward computations, we have
\begin{eqnarray*}
\Dorfman{X+\alpha,\Dorfman{Y+\beta,Z+\gamma}}&=&[X,[Y,Z]]+L_XL_Y\gamma-L_XL_Z\beta+L_X\dM
i_Z\beta-L_{[Y,Z]}\alpha+\dM i_{[Y,Z]}\alpha,\\
\Dorfman{\Dorfman{X+\alpha,Y+\beta},Z+\gamma}&=&[[X,Y],Z]+L_{[X,Y]}\gamma-L_ZL_X\beta+L_ZL_Y\alpha-L_Z\dM
i_Y\alpha\\
&&+\dM i_ZL_X\beta-\dM i_ZL_Y\alpha+\dM i_Z\dM i_Y\alpha,\\
\Dorfman{Y+\beta,\Dorfman{X+\alpha,Z+\gamma}}&=&[Y,[X,Z]]+L_YL_X\gamma-L_YL_Z\alpha+L_Y\dM
i_Z\alpha-L_{[X,Z]}\beta+\dM i_{[X,Z]}\beta.
 \end{eqnarray*}
Thus we only need to show that
$$
L_X\dM i_Z\beta+\dM i_{[Y,Z]}\alpha=L_Y\dM i_Z\alpha+\dM
i_{[X,Z]}\beta-L_Z\dM i_Y\alpha+\dM i_ZL_X\beta-\dM i_ZL_Y\alpha+\dM
i_Z\dM i_Y\alpha,
$$
which is equivalent to that
\begin{eqnarray*}
L_X\dM i_Z\beta&=&\dM i_{[X,Z]}\beta+\dM i_ZL_X\beta,\\
\dM i_{[Y,Z]}\alpha&=&L_Y\dM i_Z\alpha+-L_Z\dM i_Y\alpha-\dM
i_ZL_Y\alpha+\dM i_Z\dM i_Y\alpha.
 \end{eqnarray*}
 It is straightforward to see that they follow from the Cartan
 formula $L_X=i_X\circ\dM+\dM\circ i_X$ and
\begin{eqnarray*}
\dM\circ L_X&=&L_X\circ\dM,\\
i_{[X,Y]}\alpha&=&i_XL_Y\alpha-L_Yi_X\alpha=-(i_YL_X\alpha-L_Xi_Y\alpha).
  \end{eqnarray*}

At last, $\rho\Dorfman{X+\alpha,Y+\beta}=[X,Y]$ is obvious.
Therefore, $(TM\oplus\wedge^nT^*M,\Dorfman{\cdot,\cdot},\rho)$ is a
Leibniz algebroid.

(3). With the same notations as above, the left hand side of
(\ref{eqn:relation}) is equal to
$$
\half L_Xi_Y\gamma+\half L_Xi_Z\beta.
$$
The right hand side is equal to
$$
\half(i_{[X,Y]}\gamma+i_ZL_X\beta+i_YL_X\gamma+i_{[X,Z]}\beta).
$$
The conclusion follows from
$$
i_{[X,Y]}=L_Xi_Y-i_YL_X.
$$
This finishes the proof. \qed
\begin{rmk}
Recently several generalized Courant algebroid structures are
introduced such as $E$-Courant algebroids \cite{clsecourant} and
Loday algebroids \cite{Loday}. $E$-Courant algebroids were
introduced to unify Courant-Jacobi algebroids and omni-Lie
algebroids \cite{clomni}.  In an $E$-Courant algebroid, the pairing
takes value in a vector bundle $E$ and the target of the anchor is
the covariant differential operator bundle $\dev E$. Since the
pairing in $\mathcal{T}^n$ is $\wedge^{n-1}T^*M$-valued, we expect
it to be a $\wedge^{n-1}T^*M$-Courant algebroid.
 Unfortunately this is not true since
the composition of the Lie derivation and $\rho$ is not a bundle
map. In a Loday algebroid, the compatibility condition between the
pairing and the bracket is removed. Since the pairing in
$\mathcal{T}^n$ takes value in $\wedge^{n-1}T^*M$, it fails to be a
Loday algebroid. It is a future work to pursue the geometry behind
this higher-order bracket.
\end{rmk}

\vspace{2mm}

The standard Courant bracket can be deformed by a closed 3-form
\cite{severa3form}. Similarly, consider the following deformed
Dorfman bracket of order $n$ by an $(n+2)$-form $\Theta$,
\begin{equation}\label{CDbracketdeformed}
\Dorfman{X+\alpha,Y+\beta}_\Theta=[X,Y]+L_X\beta-L_Y\alpha+\dM
i_Y\alpha+i_{X\wedge Y}\Theta.
\end{equation}
Similar as the classical case, we have
\begin{pro}
The deformed bracket $\Dorfman{\cdot,\cdot}_\Theta$ satisfies the
Leibniz rule iff $\Theta$ is closed.
\end{pro}

For any $\Phi\in\Omega^{n+1}(M)$, define
$e^\Phi:\mathcal{T}^n\longrightarrow \mathcal{T}^n$ by
$$
e^\Phi(X+\alpha)=X+\alpha+i_X\Phi.
$$
It is straightforward to obtain
\begin{pro}\label{pro:transformation}
For any $X+\alpha,Y+\beta\in\Gamma(\mathcal{T}^n)$, we have
$$
e^\Phi\Dorfman{X+\alpha,Y+\beta}_{\dM\Phi}=\Dorfman{e^\Phi(X+\alpha),e^\Phi(Y+\beta)}.
$$
Consequently,  $e^\Phi$ is an isomorphism from the Leibniz algebroid
$(\mathcal{T}^n,\Dorfman{\cdot,\cdot}_{\dM \Phi},\rho)$ to the
Leibniz algebroid $(\mathcal{T}^n,\Dorfman{\cdot,\cdot},\rho)$.
\end{pro}

Furthermore, it is obvious that for any differeomorphism
$f\in\rm{Diff}(M)$, $(f_*,{f^*}^{-1}):\mathcal{T}^n\longrightarrow
\mathcal{T}^n$ is an automorphism preserving the higher-order
Dorfman bracket $\Dorfman{\cdot,\cdot}$. By Proposition
\ref{pro:transformation}, if an $(n+1)$-form $\Phi$ is closed,
$e^\Phi$ is also an automorphism preserving the higher-order Dorfman
bracket $\Dorfman{\cdot,\cdot}$. Similar as the classical case, we
have

\begin{thm}The group of orthogonal automorphisms of $\mathcal{T}^n$ preserving
the higher-order Dorfman bracket $\Dorfman{\cdot,\cdot}$ is the
semidirect product of $\rm{Diff}(M)$ and $\Omega^{n+1}_{\rm{c}}(M)$,
where $\Omega^{n+1}_{\rm{c}}(M)$ is the set of closed $(n+1)$-forms.
\end{thm}
We omit the proof which is similar as the proof of Proposition 3.24
in \cite{gualtieri}.\vspace{3mm}

In the classical case, Dirac structures of the Courant algebroid
$TM\oplus T^*M$ unify both presymplectic structures and Poisson
structures. More precisely, the graph of  a 2-form
$\omega\in\Omega^2(M)$ is a Dirac structure iff $\omega$ is closed
and the graph of a 2-vector field $\pi\in\frkX^2(M)$ is a Dirac
structure iff $\pi$ is a Poisson structure. Here a Dirac structure
means a maximal isotropic subbundle whose section space is closed
under the Courant bracket. See \cite{Co} for more details.

In the higher-order case, similarly we have that the graph of an
$(n+1)$-form is closed under the higher-order Courant bracket iff
$\omega$ is a premultisymplectic structure, i.e. $\dM w=0$. The
graph of an $(n+1)$-vector field is closed under the higher-order
Courant bracket iff $\pi$ is a Nambu-Poisson structure. However, the
graph of  an $(n+1)$-vector field is not isotropic, thus we can not
talk about Dirac structures. In  \cite{NambuDirac}, the author
introduced Nambu Dirac structures to unify both the structures. In
the following two sections, we will study Nambu-Poisson structures
and multisymplectic structures respectively.

\section{Nambu-Poisson structures}
A Poisson structure on an $m$-dimensional smooth manifold $M$ is a
bilinear skew-symmetric map
$\{\cdot,\cdot\}:\CWM\times\CWM\longrightarrow\CWM$ such that the
Leibniz rule and the Jacobi identity are satisfied. It is well known
that it is equivalent to a bivector field $\pi\in\wedge^2\frkX(M)$
such that  $[\pi,\pi]=0$. The relation is given by $\{f,g\}=\pi(\dM
f,\dM g)$. We usually denote a Poisson manifold by $(M,\pi)$. For
any $\sigma\in\frkX^{n+1}(M))$ and $\theta\in\Omega^{n+1}(M))$,
$\sigma^\sharp:\wedge^nT^*M\longrightarrow TM$ and
$\theta_\sharp:TM\longrightarrow \wedge^nT^*M$ are given by
$$
\sigma^\sharp(\xi)=i_\xi\sigma,\quad\theta_\sharp(X)=i_X\theta,\quad\forall~\xi\in\Omega^n(M),\forall~X\in\frkX(M).
$$
Denote by $\frkG_\sigma\subset TM\oplus\wedge^nT^*M$ (resp.
$\frkG_\theta$) the graph of $\sigma^\sharp$ (resp.
$\theta_\sharp$).

For a Poisson manifold $(M,\pi)$, there is a skew-symmetric bracket
operation $[\cdot,\cdot]_\pi$ on $1$-forms which is given by
\begin{equation}\label{defi:bracketpi}
[\alpha,\beta]_\pi=L_{\pi^\sharp(\alpha)}\beta-L_{\pi^\sharp(\beta)}\alpha+\dM
i_{\pi^\sharp(\beta)}\alpha,\quad\forall~\alpha,\beta\in\Omega^1(M).
\end{equation}
It follows that $(T^*M,[\cdot,\cdot]_\pi,\pi^\sharp)$ is a Lie
algebroid. Also we have $ [\dM f,\dM g]_\pi=\dM\{f,g\}. $

\begin{defi}{\rm\cite{TakhtajanNambu}}
A Nambu-Poisson structure of order $n$ on $M$ is an $(n+1)$-linear
map
$\{\cdot,\cdots,\cdot\}:\CWM\times\cdots\times\CWM\longrightarrow\CWM$
satisfying the following properties:
\begin{itemize}
\item[\rm(1)] skewsymmetry, i.e.
$$
\{f_1,\cdots,f_{n+1}\}=(-1)^{|\sigma|}\{f_{\sigma(1)},\cdots,f_{\sigma(n+1)}\},
$$
for all $f_1,\cdots,f_{n+1}\in \CWM$ and $\sigma\in
\mathrm{Symm}(n+1);$
\item[\rm(2)]  the Leibniz rule,
i.e. for any $g\in\CWM$, we have
$$
\{f_1g,\cdots,f_{n+1}\}=f_1\{g,\cdots,f_{n+1}\}+g\{f_1,\cdots,f_{n+1}\};
$$

\item[\rm(3)] integrability condition:
$$\{f_1,\cdots,f_{n},\{g_1,\cdots,g_{n+1}\}\}=\sum_{i=1}^{n+1}\{g_1,\cdots,\{f_1,\cdots,f_n,g_i\},\cdots,g_{n+1}\}.$$
\end{itemize}
\end{defi}

In particular, a Nambu-Poisson structure of order $1$ is exactly a
usual Poisson structure. Similar as the fact that a Poisson
structure is determined by a bivector field, a Nambu-Poisson
structure of order $n$ is determined by an $(n+1)$-vector field
$\pi\in\frkX^{n+1}(M)$ such that
$$
\{f_1,\cdots,f_{n+1}\}=\pi(\dM f_1,\cdots,\dM f_{n+1}).
$$
An $(n+1)$-vector field $\pi\in\frkX^{n+1}(M)$ is a Nambu-Poisson
structure if and only if for $f_1,\cdots,f_n\in\CWM$, we have
$$
L_{\pi^\sharp(\dM f_1\wedge\cdots\wedge\dM f_n)}\pi=0.
$$

\begin{lem}
Let $\pi\in\frkX^{n+1}(M)$ be a Nambu-Poisson structure, the
following two statements are equivalent:
\begin{itemize}
\item[\rm(i)]The graph of $\pi^\sharp$ is closed under the higher-order
Courant bracket operation \eqref{Cbracket};
\item[\rm(ii)]The graph of $\pi^\sharp$ is closed under the higher-order
Dorfman bracket operation \eqref{CDbracket}.
\end{itemize}
\end{lem}
\pf By (\ref{CCD}), we only need to prove that
$$
\pi^\sharp(d(i_{\pi^\sharp(\alpha)}\beta+i_{\pi^\sharp(\beta)}\alpha))=\pi^\sharp(L_{\pi^\sharp(\alpha)}\beta-i_{\pi^\sharp(\alpha)}d\beta
+L_{\pi^\sharp(\beta)}\alpha-i_{\pi^\sharp(\beta)}d\alpha)=0,
$$
for any $\alpha,\beta\in\Omega^n(M)$. The conclusion follows from
the following two formulas:
\begin{eqnarray*}
  \pi^\sharp(L_{\pi^\sharp(\alpha)}\beta)&=&[\pi^\sharp(\alpha),\pi^\sharp(\beta)]+(-1)^n(i_{d\alpha}\pi)\pi^\sharp(\beta);\\
\pi^\sharp(i_{\pi^\sharp(\alpha)}d\beta)&=&(-1)^n(i_{d\beta}\pi)\pi^\sharp(\alpha).
\end{eqnarray*}
The proof is finished. \qed

\begin{thm}
Let $\pi$ be an $(n+1)$-vector field, the graph
$\frkG_\pi\subset\mathcal{T}^n$  is closed under the higher-order
Dorfman bracket $\Dorfman{\cdot,\cdot}$ iff $\pi$ is a Nambu-Poisson
structure.
\end{thm}
\pf We have
\begin{eqnarray*}
\Dorfman{\pi^\sharp(\alpha)+\alpha,\pi^\sharp(\beta)+\beta}=[\pi^\sharp(\alpha),\pi^\sharp(\beta)]+L_{\pi^\sharp(\alpha)}\beta-L_{\pi^\sharp(\beta)}\alpha+\dM
i_{\pi^\sharp(\beta)}\alpha.
\end{eqnarray*}
Thus the graph $\frkG_\pi$ is closed iff
\begin{eqnarray*}
~[\pi^\sharp(\alpha),\pi^\sharp(\beta)]&=&\pi^\sharp(L_{\pi^\sharp(\alpha)}\beta-L_{\pi^\sharp(\beta)}\alpha+\dM
i_{\pi^\sharp(\beta)}\alpha)\\
&=&\pi^\sharp(L_{\pi^\sharp(\alpha)}\beta-i_{\pi^\sharp(\beta)}\dM
\alpha).
\end{eqnarray*}
On the other hand, we have
\begin{eqnarray*}
~[\pi^\sharp(\alpha),\pi^\sharp(\beta)]&=&L_{\pi^\sharp(\alpha)}\langle\beta,\pi\rangle\\
&=&\langle
L_{\pi^\sharp(\alpha)}\beta,\pi\rangle+\langle\beta,L_{\pi^\sharp(\alpha)}\pi\rangle.
\end{eqnarray*}
Therefore, by the equality
$$\pi^\sharp(L_{\pi^\sharp(\alpha)}\beta-i_{\pi^\sharp(\beta)}\dM
\alpha)=[\pi^\sharp(\alpha),\pi^\sharp(\beta)],$$ we have
\begin{eqnarray*}
\langle
L_{\pi^\sharp(\alpha)}\pi,\beta\rangle=-\pi^\sharp(i_{\pi^\sharp(\beta)}\dM\alpha).
\end{eqnarray*}
Let $\alpha=\dM f_1\wedge\cdots\wedge\dM f_n$, then $\dM \alpha=0$.
Thus for any $n$-form $\beta$, we have
$$
\langle L_{\pi^\sharp(\dM f_1\wedge\cdots\wedge\dM
f_n)}\pi,\beta\rangle=0,
$$
which implies that
$$
 L_{\pi^\sharp(\dM f_1\wedge\cdots\wedge\dM
f_n)}\pi=0,
$$
i.e. $\pi$ is a Nambu-Poisson structure of order $n$.

Conversely, if $\pi$ is a Nambu-Poisson structure of order $n$, by
the local expression of
 $\pi$ \cite{MarrerodynamicsNambu,MarmonPoisson,NakanishiNM}, at any point $x\in M$
satisfying $\pi(x)\neq0$, we can choose a local coordinates
$(x_1,\cdots,x_{n+1},\cdots,x_m)$ around $x$ such that
$$
\pi=\frac{\partial}{\partial
x_1}\wedge\cdots\wedge\frac{\partial}{\partial x_{n+1}}.
$$
Then it is not hard to see that the graph $\frkG_\pi$  is closed
under the higher-order Dorfman bracket. \qed

\begin{cor}
Let $\pi$ be a Nambu-Poisson structure of order $n$.
\begin{itemize}
\item[\rm(1).] The graph  $\frkG_\pi\subset \mathcal{T}^n$  is a sub-Leibniz algebroid
of the Leibniz algebroid
$(\mathcal{T}^n,\Dorfman{\cdot,\cdot},\rho)$.
\item[\rm(2).] There is an induced bracket operation $[\cdot,\cdot]_\pi$ on $
\Omega^{n}(M)$,
\begin{equation}\label{defi:bracketpin}
[\alpha,\beta]_\pi=L_{\pi^\sharp(\alpha)}\beta-L_{\pi^\sharp(\beta)}\alpha+\dM
i_{\pi^\sharp(\beta)}\alpha,
\end{equation}
which is a natural generalization of the bracket given by
\eqref{defi:bracketpi}. Consequently,
$(\wedge^nT^*M,[\cdot,\cdot]_\pi,\pi^\sharp)$ is a Leibniz
algebroid.
\end{itemize}
\end{cor}

   In the following, we prove that there is also a Leibniz algebra
$(\Omega^{n-1}(M),\{\cdot,\cdot\}_\pi)$ associated to any
Nambu-Poisson structure $\pi$ of order $n$.

For any $\xi,\eta\in\Omega^{n-1}(M)$, we have
$$
~[\dM \xi,\dM\eta]_\pi=L_{\pi^\sharp\dM\xi}\dM\eta=\dM
L_{\pi^\sharp(\dM\xi)}\eta.
$$
Define the bilinear bracket $\{\cdot,\cdot\}_\pi$ on
$\Omega^{n-1}(M)$ by
\begin{equation}
\{\xi,\eta\}_\pi=L_{\pi^\sharp(\dM\xi)}\eta.
\end{equation}
Obviously we have
\begin{equation}\label{eqn:equal}
[\dM \xi,\dM\eta]_\pi=\dM\{\xi,\eta\}_\pi.
\end{equation}
\begin{thm}
Let $\pi$  be a Nambu-Poisson structure  of order $n$. Then
$(\Omega^{n-1}(M),\{\cdot,\cdot\}_\pi)$ is a Leibniz algebra. In
particular, if $n=1$, $(\CWM,\{\cdot,\cdot\}_\pi)$ is the canonical
Poisson algebra.
\end{thm}
\pf First by (\ref{eqn:equal}), we have
$$
\pi^\sharp \dM
L_{\pi^\sharp(\dM\xi)}\eta=\pi^\sharp\dM\{\xi,\eta\}_\pi=\pi^\sharp[\dM
\xi,\dM\eta]_\pi.
$$
The fact that $(\wedge^{n}T^*M,[\cdot,\cdot]_\pi,\pi^\sharp)$ is a
Leibniz algebroid yields that
$$
\pi^\sharp \dM L_{\pi^\sharp(\dM\xi)}\eta=[\pi^\sharp(\dM
\xi),\pi^\sharp(\dM\eta)].
$$
Thus we have
\begin{eqnarray*}
&&\{\xi,\{\eta,\gamma\}_\pi\}_\pi-\{\{\xi,\eta\}_\pi,\gamma\}_\pi-\{\eta\{\xi,\gamma\}_\pi\}_\pi\\&=&L_{\pi^\sharp(\dM\xi)}L_{\pi^\sharp(\dM\eta)}\gamma
-L_{\pi^\sharp(\dM L_{\pi^\sharp(d\xi)}\eta)}\gamma-L_{\pi^\sharp(\dM\eta)}L_{\pi^\sharp(\dM\xi)}\gamma\\
&=&L_{[\pi^\sharp(\dM\xi),\pi^\sharp(\dM\eta)]}\gamma-L_{\pi^\sharp(\dM
L_{\pi^\sharp(d\xi)}\eta)}\gamma\\
&=&0,
\end{eqnarray*}
which implies that $(\Omega^{n-1}(M),\{\cdot,\cdot\}_\pi)$ is a
Leibniz algebra. The other conclusion is obvious and the proof is
completed. \qed

\begin{rmk}  In
\cite{MarreroLeibnizNambu}, the authors introduced a  bracket
operation $[\cdot,\cdot]^\pi:\Omega^{n}(M)\times
\Omega^{n}(M)\longrightarrow \Omega^{n}(M)$ associated with  a
Nambu-Poisson structure $\pi$ of order $n$ by setting
\begin{equation}\label{defi:bracketleibniz}
[\alpha,\beta]^\pi=L_{\pi^\sharp(\alpha)}\beta+(-1)^{n+1}(i_{\dM\alpha}\pi)\beta.
\end{equation}
This bracket is not skew-symmetric. However, the Leibniz rule is
satisfied, i.e.
$$
[\alpha,[\beta,\gamma]^\pi]^\pi=[[\alpha,\beta]^\pi,\gamma]^\pi+[\beta,[\alpha,\gamma]^\pi]^\pi.
$$
Furthermore, we have
$$
\pi^\sharp[\alpha,\beta]^\pi=[\pi^\sharp(\alpha),\pi^\sharp(\beta)].
$$
It follows that $(\wedge^{n}T^*M,[\cdot,\cdot]^\pi,\pi^\sharp)$ is a
Leibniz algebroid.

It is obvious that the bracket $[\cdot,\cdot]_\pi$ given by
\eqref{defi:bracketpin} and the bracket $[\cdot,\cdot]^\pi$ given by
\eqref{defi:bracketleibniz} are not the same even though both of
them define Leibniz algebroid structures on $\wedge^nT^*M$. By the
equalities
$$
\pi^\sharp[\alpha,\beta]^\pi=[\pi^\sharp(\alpha),\pi^\sharp(\beta)]\quad
\pi^\sharp[\alpha,\beta]_\pi=[\pi^\sharp(\alpha),\pi^\sharp(\beta)],
$$
we have$$ \pi^\sharp([\alpha,\beta]_\pi-[\alpha,\beta]^\pi)=0.
$$
However, we believe that our bracket \eqref{defi:bracketpin} is much
more reasonable  since it is induced from the higher-order Courant
bracket \eqref{CDbracket} and  it is also a generalization of
\eqref{defi:bracketpi}.
\end{rmk}

\section{Multisymplectic structures }

\begin{defi}\label{defi:2-plectic}
An (n+1)-form $\omega$ on a smooth manifold $M$ is a multisymplectic
structure of order $n$, or an $n$-plectic structure if it is closed
($\dM \omega=0$) and nondegenerate:
$$
i_X\omega=0\Rightarrow X=0,\quad\forall~X\in T_mM,~m\in M.
$$
If $\omega$ is an n-plectic structure on $M$, we call the pair
$(M,\omega)$ an n-plectic manifold.
\end{defi}

Remove the nondegenerate condition, we get premultisymplectic
structures. About more details about multisymplectic structures,
please see \cite{baez:classicalstring} and the references there.

\begin{thm}
Let $\omega$ be an $(n+1)$-form, the graph $\frkG_\omega $ is a
maximal isotropic subbundle of $\mathcal{T}^n$. Furthermore, it is
closed under the higher-order Dorfman bracket \eqref{CDbracket} iff
$\omega$ is closed.
\end{thm}
\pf  For any $n$-plectic structure $\omega$, we have
$$
\pairing{X+\omega_\sharp(X),Y+\omega_\sharp(Y)}=\half
(i_Yi_X\omega+i_Xi_Y\omega)=0.
$$
which implies that $\frkG_\omega\subset TM\oplus\wedge^nT^*M$ is
isotropic.

On the other hand, since $\omega$ is  closed, we have
\begin{equation}\label{eqn:closed}
i_{[X,Y]}\omega=L_Xi_Y\omega- i_Y\dM i_X\omega,
\end{equation} which yields that $\frkG_\omega$ is closed under the
bracket operation $\Dorfman{\cdot,\cdot}$. \qed\vspace{3mm}

Consider the deformed bracket (\ref{CDbracketdeformed}), we have

\begin{pro}
For an $(n+1)$-form $\omega$, its graph $\frkG_\omega\subset
\mathcal{T}^n$ is closed under the deformed bracket
\eqref{CDbracketdeformed} if and only if $$\dM \omega+\Theta=0.$$
\end{pro}
\pf It is straightforward to see that the graph of $\omega$ is
closed under the deformed bracket operation
$\Dorfman{\cdot,\cdot}_\Theta$ is equivalent to that
$$
i_{[X,Y]}\omega=L_Xi_Y\omega-i_Y\dM i_X\omega+\Theta(X,Y).
$$
On the other hand, we have
$$
i_{[X,Y]}\omega=L_Xi_Y\omega-i_YL_X\omega=L_Xi_Y\omega-i_Yi_X\dM
\omega-i_Y\dM i_X\omega.
$$
Thus, the graph of $\omega$ is closed under the deformed bracket
(\ref{CDbracketdeformed}) if and only if
$$
i_Yi_X(\dM \omega+\Theta)=0,\quad \forall~X,Y\in\frkX(M),
$$
which is equivalent to that $\dM \omega+\Theta=0$. \qed\vspace{3mm}

By Proposition \ref{pro:transformation}, we have
\begin{cor}
If there is an $(n+1)$-form $\omega$ such that its graph
$\frkG_\omega\subset \mathcal{T}^n$ is closed under the deformed
bracket \eqref{CDbracketdeformed}, then $e^{-\omega}$ is an
isomorphism from the Leibniz algebroid
$(\mathcal{T}^n,\Dorfman{\cdot,\cdot}_\Theta,\rho)$ to the Leibniz
algebroid $(\mathcal{T}^n,\Dorfman{\cdot,\cdot},\rho)$
\end{cor}

 For an $n$-plectic manifold $(M,\omega)$, an $(n-1)$-form
$\xi\in\Omega^{n-1}(M)$ is Hamiltonian if there is a vector field
$X_\xi$ such that
\begin{equation}\label{def:H}
\dM \xi=i_{X_\xi}\omega.
\end{equation} We  say $X_\xi$ is the
Hamiltonian vector field corresponding to $\xi$. The set of
Hamiltonian $(n-1)$-forms and the set of Hamiltonian vector fields
on $(M,\omega)$ are denoted as $\Ham(M)$ and $\frkX_H(M)$.

In \cite{baez:classicalstring}, Baez and Rogers generalize the
Poisson bracket of functions in symplectic geometry to a bracket of
Hamiltonian $(n-1)$-forms in two ways. One is the so called
hemi-bracket $\{\cdot,\cdot\}_h$ defined as follows:
\begin{equation}\label{defi:brackethemi}
\{\xi,\eta\}_h=L_{X_\xi}\eta,\quad\xi,\eta\in\Ham(M).
\end{equation}
The other one is the so called semi-bracket $\{\cdot,\cdot\}_s$
defined as follows:
$$
\{\xi,\eta\}_s=i_{X_\xi}i_{X_\eta}\omega.
$$
However, both of them are not Lie brackets. In particular, for
2-plectic structures,  the author obtained  hemistrict Lie
2-algebras and  semistrict Lie 2-algebras.

Lie 2-algebras are categorification of Lie algebras. The Jacobi
identity is replaced by a natural isomorphism, called the
Jacobiator, which also satisfies a certain law of its own. A good
introduction for this subject is \cite{baez:2algebras}.




\begin{rmk}
Note that by \eqref{def:H}, usually $f\xi$ is not a Hamiltonian
$(n-1)$-form even if $\xi$ is a Hamiltonian $(n-1)$-form. Thus,
$\Ham(M)$ is not the section space of a vector bundle.
\end{rmk}

In symplectic geometry and Poisson geometry, except the Poisson
bracket of functions,  one can also define a bracket on $1$-forms
and obtain the cotangent bundle Lie algebroid. In the following, we
generalize this idea to the case of $n$-plectic geometry. For an
$n$-plectic manifold $(M,\omega)$, not all the $(n-1)$-forms are
Hamiltonian $(n-1)$-forms. So there is no way to expect that we can
introduce a bracket on all the $n$-forms such that $\wedge^nT^*M$ is
a Lie algebroid. However, this can be done for a subbundle of
$\wedge^nT^*M$.

\begin{defi}
Let $(M,\omega)$ be an $n$-plectic manifold. $\alpha\in\Omega^n(M)$
is called an admissible $n$-form if there is a vector field
$X_\alpha\in\frkX(M)$ such that $ \alpha=i_{X_\alpha}\omega.
$
\end{defi}

It is straightforward to see that $\xi$ is a Hamiltonian
$(n-1)$-form iff $\dM\xi$ is an admissible $n$-form. Obviously, if
$\alpha$ is an admissible $n$-form, then $f\alpha$ is also an
admissible $n$-form since we have
$$f\alpha=fi_{X_\alpha}\omega=i_{fX_\alpha}\omega.$$
Therefore, the set of admissible $n$-forms is the  section space of
a subvector bundle of $\wedge^nT^*M$, which we call the {\bf
admissible bundle} and denote by $A$. Since $\omega $ is
nondegenerate, we obtain a bundle map
$\omega^\natural:A\longrightarrow TM$ by
\begin{equation}\label{def:anchoronA}
\omega^\natural(\alpha)=X_\alpha\quad \alpha\in \Gamma(A).
\end{equation}

Define a skew-symmetric bracket operation $[\cdot,\cdot]_\omega$ on
$\Gamma(A)$ by setting
\begin{equation}\label{def:bracketonA}
[\alpha,\beta]_\omega=L_{X_\alpha}\beta-L_{X_\beta}\alpha-\dM
i_{X_\alpha}i_{X_\beta}\omega.
\end{equation}
It is not hard to deduce that
\begin{lem}\label{lem:bracket2}
With the above notations, for any admissible $n$-forms
$\alpha,\beta$, their bracket $[\alpha,\beta]_\omega$ is again an
admissible $n$-form. Moreover, we have
$$
X_{[\alpha,\beta]_\omega}=[X_\alpha,X_\beta].
$$
Thus, the bracket \eqref{def:bracketonA} is well defined.
\end{lem}

\begin{thm}
With the above notations, $(A,[\cdot,\cdot]_\omega,\omega^\natural)$
is a Lie algebroid with the anchor $\omega^\natural$, where
$[\cdot,\cdot]_\omega$ and $\omega^\natural$ are given by
\eqref{def:bracketonA} and \eqref{def:anchoronA} respectively.
\end{thm}
\pf By Lemma \ref{lem:bracket2}, we have
$$
[[\alpha,\beta]_\omega,\gamma]_\omega+c.p.=i_{[[X_\alpha,X_\beta],X_\gamma]+c.p.}\omega=0,
$$
which implies that $[\cdot,\cdot]_\omega$ is a Lie bracket.
Furthermore, for any $f\in \CWM$, we have
\begin{eqnarray*}
[\alpha,f\beta]_\omega&=&i_{[X_\alpha,fX_\beta]}\omega=fi_{[X_\alpha,X_\beta]}\omega+X_\alpha(f)i_{X_\beta}\omega\\
&=&f[\alpha,\beta]_\omega+X_\alpha(f)\beta\\
&=&f[\alpha,\beta]_\omega+\omega^\natural(\alpha)(f)\beta,
\end{eqnarray*}
which implies that $(A,[\cdot,\cdot]_\omega,\omega^\natural)$ is a
Lie algebroid.\qed\vspace{3mm}
\begin{rmk}
When $n=1$, $A$ is exactly the cotangent bundle $T^*M$ and we
recover the cotangent bundle Lie algebroid.
\end{rmk}

In the following, we consider the bracket of exact admissible
$n$-forms. We will see how  we obtain the hemi-bracket
(\ref{defi:brackethemi}) on $\Ham(M)$  from our Lie algebroid
$(A,[\cdot,\cdot]_\omega,\omega^\sharp)$.

Let $\alpha=\dM \xi,~\beta=\dM \eta$ be two exact admissible
$n$-forms, where $\xi,\eta$ are $(n-1)$-forms. Then we have
\begin{equation}\label{eqn:exact}
[\alpha,\beta]_\omega=[\dM \xi,\dM \eta]_\omega=L_{X_{\dM\xi}}\dM
\eta=\dM L_{X_{\dM\xi}} \eta.
\end{equation} Thus, $[\dM \xi,\dM \eta]_\omega$ is
again an exact admissible $n$-form.
   In this case we simply denote
$X_{\dM\xi}$ by $X_\xi$. It is also the same as the Hamiltonian
vector field of the Hamiltonian $(n-1)$-form $\xi$. By
(\ref{eqn:exact}), there is an induced bracket $\{ \cdot,
\cdot\}_\omega$ of Hamiltonian $(n-1)$-forms by setting
$$
\{ \xi, \eta\}_\omega=L_{X_\xi}\eta.
$$
This is exactly the hemi-bracket (\ref{defi:brackethemi}) of
Hamiltonian $(n-1)$-forms. Obviously, we have
$$\dM \{ \xi, \eta\}_\omega =[\dM \xi,\dM \eta]_\omega.$$ However,
even though the bracket operation $[\cdot,\cdot]_\omega$ is
skew-symmetric,  $\{ \cdot, \cdot\}_\omega$ may be not
skew-symmetric. The obstruction is given by an exact term,  see
\cite{baez:classicalstring} for more details. In particular, for
2-plectic structures,  we obtain a hemistrict Lie 2-algebra of which
the degree-0 part is $\Ham(M)$ and the degree-1 part is $\CWM$.

\begin{rmk} Rogers proved that one can associate a Lie $n$-algebra to
an $n$-plectic structures \cite{Rogers}. It will be interesting to
study whether one can associate $L_\infty$-algebras to Nambu-Poisson
structures.
\end{rmk}


\begin{thebibliography}{99}

\bibitem{alekseevmanin}
Alekseev A and Kosmann-Schwarzbach Y.
\newblock Manin pairs and moment maps.
\newblock { J. Differential Geom.}, 56(1):133--165, 2000.

\bibitem{baez:2algebras}
Baez J C and  Crans A S.
\newblock Higher-dimensional algebra. {VI}. {L}ie 2-algebras.
\newblock { Theory Appl. Categ.}, 12:492--538 (electronic), 2004

\bibitem{baez:classicalstring}
Baez J C,  Hoffnung A E, and  Rogers C L.
\newblock Categorified symplectic geometry and the classical string.
\newblock { Comm. Math. Phys.}, 293(3):701--725, 2010

\bibitem{chatterjee2}
Chatterjee R.
\newblock Dynamical symmetries and {N}ambu mechanics.
\newblock { Lett. Math. Phys.}, 36(2):117--126, 1996

\bibitem{chatterjee1}
Chatterjee R and Takhtajan L.
\newblock Aspects of classical and quantum {N}ambu mechanics.
\newblock { Lett. Math. Phys.}, 37(4):475--482, 1996

\bibitem{clomni}
Chen Z and Liu Z J.
\newblock Omni-{L}ie algebroids.
\newblock { J. Geom. Phys.}, 60(5):799--808, 2010

\bibitem{clsecourant}
Chen Z,  Liu Z J, and Sheng Y H.
\newblock {$E$-Courant algebroids}, Int. Math. Res. Not. Vol. (2010), No. 22, pp. 4334-4376.

\bibitem{Co}
Courant T J.
\newblock Dirac manifolds.
\newblock { Trans. Amer. Math. Soc.}, 319(2):631--661, 1990

\bibitem{GrabowskiNambu2}
Grabowski J and Marmo G.
\newblock On {F}ilippov algebroids and multiplicative {N}ambu-{P}oisson
  structures.
\newblock { Differential Geom. Appl.}, 12(1):35--50, 2000

\bibitem{gualtieri}
Gualtieri M.
\newblock {Generalized complex geometry}, arxiv:math.DG/0401221

\bibitem{NambuDirac}
Hagiwara Y.
\newblock Nambu-{D}irac manifolds.
\newblock { J. Phys. A}, 35(5):1263--1281, 2002

\bibitem{MarrerodynamicsNambu}
Ib{\'a}{\~n}ez R, de~Le{\'o}n M, Marrero J, and Mart{\'{\i}}n~de
  Diego D.
\newblock Dynamics of generalized {P}oisson and {N}ambu-{P}oisson brackets.
\newblock {J. Math. Phys.}, 38(5):2332--2344, 1997

\bibitem{MarreroLeibnizNambu}
Ib{\'a}{\~n}ez R, de~Le{\'o}n M,  Marrero J, and Padr{\'o}n E.
\newblock Leibniz algebroid associated with a {N}ambu-{P}oisson structure.
\newblock { J. Phys. A}, 32(46):8129--8144, 1999

\bibitem{Liblandmeinrenken}
Li-Bland D and Meinrenken E.
\newblock Courant algebroids and {P}oisson geometry.
\newblock { Int. Math. Res. Not. IMRN}, (11):2106--2145, 2009

\bibitem{lwx}
Liu Z J, Weinstein A, and Xu P.
\newblock Manin triples for {L}ie bialgebroids.
\newblock { J. Diff. Geom.}, 45(3):547--574, 1997

\bibitem{MarmonPoisson}
Marmo G, Vilasi G, and  Vinogradov A M.
\newblock The local structure of {$n$}-{P}oisson and {$n$}-{J}acobi manifolds.
\newblock { J. Geom. Phys.}, 25(1-2):141--182, 1998

\bibitem{NakanishiNM}
Nakanishi N.
\newblock On {N}ambu-{P}oisson manifolds.
\newblock { Rev. Math. Phys.}, 10(4):499--510, 1998

\bibitem{Rogers}
Rogers C.
\newblock {L-infinity algebras from multisymplectic geometry}, arXiv:1005.2230

\bibitem{Roytenbergphdthesis}
Roytenberg D.
\newblock {Courant algebroids, derived brackets and even symplectic
  supermanifolds}, PhD thesis, UC Berkeley, 1999, arXiv:math.DG/9910078

\bibitem{rw}
Roytenberg D and Weinstein A.
\newblock Courant algebroids and strongly homotopy {L}ie algebras.
\newblock { Lett. Math. Phys.}, 46(1):81--93, 1998

\bibitem{severa3form}
{\v{S}}evera P and Weinstein A.
\newblock Poisson geometry with a 3-form background.
\newblock { Progr. Theoret. Phys. Suppl.}, (144):145--154, 2001.
\newblock Noncommutative geometry and string theory (Yokohama, 2001)

\bibitem{shengJacobi}

Sheng Y H. Jacobi quasi-Nijenhuis algebroids. Rep. Math. Phys. 65,
271-287, 2010

\bibitem{Loday}
Sti{\'e}non M and Xu P.
\newblock Modular classes of {L}oday algebroids.
\newblock { C. R. Math. Acad. Sci. Paris}, 346(3-4):193--198, 2008

\bibitem{xureductionofgeneralized}
Sti{\'e}non M and Xu P.
\newblock Reduction of generalized complex structures.
\newblock { J. Geom. Phys.}, 58(1):105--121, 2008

\bibitem{TakhtajanNambu}
Takhtajan L.
\newblock On foundation of the generalized {N}ambu mechanics.
\newblock { Comm. Math. Phys.}, 160(2):295--315, 1994

\bibitem{zambonCourant}
Zambon M.
\newblock {L-infinity algebras and higher analogues of Dirac structures},
  arXiv:1003.1004

\end{thebibliography}
\end{document}